\documentclass{amsart}
\usepackage{amssymb,amsthm,amsmath}
\usepackage{latexsym,mathrsfs,comment,float}
\usepackage[all]{xy}
\usepackage{seqsplit}
\usepackage{graphicx}
\usepackage[]{algorithm2e}
\usepackage[utf8]{inputenc}    % utf8 support       %!!!!!!!!!!!!!!!!!!!!
\usepackage[T1]{fontenc}

\theoremstyle{definition}
\newtheorem{definition}{{\textbf{Definition}}}[section]
\newtheorem*{thm*}{{\textbf{Theorem}}}
\newtheorem{thm}[definition]{{\textbf{Theorem}}}
\newtheorem*{definition*}{{\textbf{Definition}}}
\newtheorem{lem}[definition]{{\textbf{Lemma}}}
\newtheorem*{lem*}{{\textbf{Lemma}}}
\newtheorem{cor}[definition]{{\textbf{Corollary}}}
\newtheorem*{cor*}{{\textbf{Corollary}}}
\newtheorem{prop}[definition]{{\textbf{Proposition}}}
\newtheorem*{prop*}{{\textbf{Proposition}}}

\newtheorem*{ex*}{{\textbf{E\footnotesize XAMPLE}}}

\def\T{\mathcal{T}}
\newcommand{\Z}{\mathbb{Z}}

\def\real{\mathbb{R}}

\def\quaternion{\mathbb{H}}
\newcommand{\F}{\mathbb{F}}

\makeatletter
\def\@wgt#1{\operatorname{wgt}(#1)}
\def\@@wgt[#1]#2{\operatorname{wgt}_{#1}(#2)}
\def\wgt{\@ifnextchar[{\@@wgt}{\@wgt}}
\makeatother
\newcommand{\cat}{\mbox{cat}}

\def\tc#1{\operatorname{tc}(#1)}
\def\TC#1{\operatorname{TC}(#1)}

\newcommand{\pr}{\mbox{pr}}

\def\ev{\operatorname{ev}}
\def\map#1{\operatorname{Map}(#1)}
\def\path#1{\operatorname{\mathcal{P}}(#1)}
\def\genus#1{\operatorname{\mathfrak{genus}}(#1)}
\def\midvert{ \ \mathstrut\vrule\hskip.1pt\ }
\def\field{\mathbb{F}}

\title{On the topological complexity of $S^3/Q_8$}
\date{}
\author{Yuya~Miyata}
\address{Graduate School of Mathmatics, Kyushu University}
\email{y-miyata@math.kyushu-u.ac.jp}
\begin{document}
\maketitle
\begin{abstract}
Topological complexity was first introduced in 2003 by Michael Farber \cite{farber2003topological} as a homotopy invariant for a connected topological space $X$, denoted by $\TC{X}$.
Although the invariant is defined in terms of elementary homotopy theory using well-known Serre path fibration, not many examples are known to be determined concretely by now.
In 2010, Iwase and Sakai showed that the topological complexity of a space is a fibrewise version of a L-S category for a fibrewise space over the space.
In this paper, we determine the topological complexity of $S^3/Q_8$ using a method produced from the fibrewise view point.
\end{abstract}
%    "位相的複雑さ"とは，弧状連結な位相空間$X$に対するホモトピー不変量であり，tc($X$)と表される\cite{farber2003topological}．
%    これは$X$のある特別な被覆に関して非常に簡潔に定義される不変量であるが，現在まで，完全に決定された空間は少ない．
%    2010年に，tc($X$)が$X$から作られるファイバーワイズ空間に対するL-Sカテゴリーに一致することが\cite{iwase2010topological}によって示された．
%    私はこのファイバーワイズL-Sカテゴリー理論を用いることで$S^3/Q_8$の位相的複雑さを決定した．

\section{Intoroduction}
In 2003, Michael Farber introduced a numerical homotopy invariant called a \textit{topological complexity} in \cite{farber2003topological}.
%\cite{farber2003topological}のなかでM.~Farberは"topological complexity"と呼ばれる不変量を定義した．
But recent years, many authors working on this subject use a `reduced' version of topological complexity which is one less than the original.
In this paper, we adopt the reduced version and related invariants introduced by Michael Farber and Mark Grant \cite{farber2008robot}.
%この節では，現在用いられている topological complexity の定義と，これに関連した不変量を紹介する\cite{farber2008robot}．

The Serre path fibration $\pi  : \map{I,B} \to B \!\times\! B$ is defined by $\pi(u) = (u(0),u(1))$.
%ファイブレーション$(ev_0,ev_1) : \mbox{Map}(I,B) \to B \!\times\! B$を$\pi$によって表すことにする．
The topological complexity is defined to estimate the number of local sections of $\pi$ whose domains cover the entire base space $B \!\times\! B$.
%topological complexityはこのファイブレーションに対して，局所切断の最小数−１として定義される非常に簡潔な不変量である．
\begin{definition}%(topological complexity)
    Let $X$ be a path-connected %topological 
    space. 
    The topological complexity of $X$, denoted by $\tc{X}$, is least integer $n$ such that there exists an open covering $U_0, \dots,U_n$ of $X \!\times\! X$ each of which is equipped with a section of $\pi$ on it.
\end{definition}

By definition, it is clear that $\tc{X}$ is nothing but the Schwarz genus of the Serre path fibration $\pi$, in other words, $\tc{X}=\genus{\pi}$, where we define the Schwarz genus as one less than the original, again, as follows.
%% 今ファイブレーション$(ev_0,ev_1) : \mbox{Map}(I,B) \to B \!\times\! B$を$\pi$によって表すと，
%またtc($X$)は，以下のように定義される Schwarz genus($\mathfrak{genus}(p)$)の$p=\pi$の場合と考えることができる．

\begin{definition}
    Let $p : E \to B$ be a fibration.
    The \textit{Schwarz genus} of $p$, denoted by $\mathfrak{genus}(p)$, is the least integer $k$ such that there is an open covering $U_0, \dots, U_k$ of $B$ each of which is equipped with a section of $p$ on it.
%, i.e, a continuous map $s_i : U_i \to E$ with $p \circ s_i = \mbox{id}_{U_i}$ for each $i = 0, \dots,k$.
\end{definition}

% Topological complexity is, by definition, closely related to the Lusternik-Schnirelmann category often abbreviated as L-S category.
Topological complexity is closely related to the Lusternik-Sch\-nirelmann category often abbreviated as L-S category.
The normalised version of L-S category of $X$ a space with base point $\ast$, denoted in this paper by $\cat{X}$, can be defined in a similar manner to the topological complexity as follows:
Let $\pi_{0} = \ev_{1}\vert{\path{X}} : \path{X} \to X$, where $\path{X}=\{u \in \map{I,X} \midvert u(0)=\ast\}$.
Then we have $\cat{X} = \genus{\pi_{0}}$.

We remark that, instead of the word `reduced', the word `normalised' is used to mean ``one less than the original'', in L-S theory.
%L-Sカテゴリー(cat($X$))は$X$のPath space fibration $\pi_0 = \mbox{ev}_1 : PX \to X$のSchwarz genus($\mathfrak{genus}(\pi_0)$)であり，
Farber has pointed out, in \cite{farber2003topological}, the relationship between topological complexity and L-S category as follows:
%現在までに L-Sカテゴリーによるtc($X$)の評価が以下のように得られている\cite{farber2003topological}:
\[ \cat(X) \le \tc{X} \le \cat(X\!\times\! X) \le 2\cat(X).\]

Farber and Grant have introduced in \cite{farber2003topological,farber2008robot,farber2007symmetric}, a computable numerical invariant called TC-weight, defined in terms of Schwarz genus as follows: Let $G$ be an abelian group.
%$H^{\ast}(X;G)$の元$u$の重み(weght)はtc($X$)の評価に用いられる不変量であり Schwarz genus を用いて定義される．
%以下にweight の主な性質を述べる\cite{farber2003topological},\cite{farber2008robot},\cite{farber2007symmetric}．

\begin{definition}
Let $p : E \to B$ be a fibration.
Then, for a non-zero element $u \in H^{\ast}(B;G)$, %the weight of $u$ with respect to $\pi$ with coefficients in $G$, 
we denote by $\wgt[p]{u;G}$ the largest integer $k \ge 0$ such that $f^{\ast}(u) = 0 \in H^{\ast}(Y;G)$ among all maps $f : Y \to B$ with $\mathfrak{genus}(f^{\ast}\pi) < k$, and we denote $\wgt[p]{0;G}=\infty$.
For the Serre path fibration $\pi : \map{I,B} \to B \!\times\! B$, we often abbreviate as $\wgt{u;G}=\wgt[\pi]{u;G}$.
\end{definition}

\begin{prop}\label{cupsum1}
%    If there exists a non-zero cohomology class $u \in H^{\ast}(X \!\times\! X)$ with $\wgt{u} \ge k$, then $\tc{X} \ge k$.
For non-zero element $u \in H^{\ast}(X \!\times\! X;G)$, we have $$\tc{X} \ge \wgt{u;G}.$$
\end{prop}

\begin{prop}\label{cupsum2}
    For $u\in H^{\ast}(X \!\times\! X;G)$, we have
    \[\wgt{u;G} \ge 1 \iff u|\Delta_X = 0 \in H^{\ast}(X;G).\]
\end{prop}
    
\begin{prop}\label{cupsum3}
Let $R$ be a ring with unit.
THen for any $u,v \in H^{\ast}(X \!\times\! X;R)$, the cup product  $uv \in H^{\ast}(X \!\times\! X;R)$ satisfies
    \[ \wgt{uv;R} \ge \wgt{u;R} + \wgt{v;R}.\]
\end{prop}

Please note that we might drop the coefficient groups from the above notations, if there are no confusions.

On the other hand in 1995, James \cite{james1995introduction} introduced a fibrewise version of a L-S category for a fibrewise space.
Using the notion established in \cite{james1995introduction}, Iwase and Sakai \cite{iwase2010topological} showed that $\tc{X}$ can be interpreted as a kind of fibrewise L-S category of a fibrewise space $E$ defined as follows:
%\cite{iwase2010topological}で岩瀬，酒井はファイバーワイズ空間に対するLS-categoryがtcと一致することを示した．
%ファイバーワイズ空間については\cite{james1995introduction}を参照していただきたい．
Let $E = (X\!\times\!X,\mbox{pr}_2,X,\Delta)$ be a fibrewise space, where $\pr_{t} : X\!\times\!X \to X$ denotes the $t$-th projection and $\Delta : X \to X\!\times\!X$ is the diagonal map.
%今，空間$X$に対してファイバーワイズ空間$E = (X\!\times\!X,\mbox{pr}_2,X,\Delta)$を定義する，
%ここで$\mbox{pr}_2$は第二成分への射影であり，$\Delta : X\to X\!\times\!X$は対角写像である．
The following theorems are the key to our main result.
%次の定理，定義は重要である．

\begin{thm}[\cite{iwase2010topological}]
    For a space $X$, we have the following equality.
    \[\tc{X} = \cat_{B}^{\ast}(E)\]
\end{thm}

\begin{thm}[\cite{iwase2010topological}]\label{inj}
    Let $E$ be a fibrewise pointed space over $X$ and $m\ge 0$. Then $\cat_{B} (E) \le m$ if and only if $\mbox{id}_E : E \to E$ has a lift to $P^{m}_{B}(\Omega_{B}(E)) \xrightarrow{e_{m}^{E}} E$ in $\T_{B}$.
\end{thm}

\begin{definition}[\cite{iwase2010topological}]
    For any $u \in H^{\ast}(E;R)$, we define
    \[ \wgt[B]{u;R} = \mbox{Min}\{ m\ge 0 | (e_{m}^{E})^{\ast} (u) \neq 0\},\]
    where $e_{m}^{E}$ denotes the fibrewise map $P^{m}_{B}(\Omega_{B}(E)) \hookrightarrow P^{\infty}_{B}(\Omega_{B}(E)) \underset{\simeq}{\xrightarrow{e_{\infty}^{E}}} E$.
\end{definition}

The main result in this paper is described as follows.
%私は$X=S^3/Q_8$のtcを$wgt_X$の計算によって得た．

\begin{thm}
    Let $S^3/Q_8$ be the quotient space of $S^3$ by the standard action of the quaternion group $Q_8 \subset \mathrm{SO}(4)$. Then we obtain %the following equality.
    \[\tc{S^3/Q_8} = \cat_{B}^{\ast}(S^3/Q_8) = 6.\]
\end{thm}

We show the theorem by using the following two lemmas which is obtained by concrete computations of weight of the element in the top dimension.
%この定理は以下の二つの補題によって示される．これらは\S 3で証明する．

\begin{lem}
    Let $z$ be a cocycle representing the generator of $H^6(X;\F_2)$, where $X = S^3/Q_8\!\times\!S^3/Q_8$.
    Then %the homomorphism induced by $e_5^{X} : P^{5}_{B}(\Omega_{B}(X)) \to P^{\infty}_{B}(\Omega_{B}(X)) \simeq X$ sends $z$ to zero in $H^6(P^{5}_{B}(\Omega_{B}(X)))$ or 
    there is a cochain $u$ satisfying $\delta u = (e_{5}^{X})^{\ast}(z)$.
\end{lem}

\begin{lem}
    Let $K$ be the space $S^{\infty}/Q_8$, we have \[P^{m}_{B}\Omega_{B}(K\times K) \simeq_{B} P^{m}_{B}\widehat{K}.\]
\end{lem}

The paper is organised as follows. In Section 2, we determine the ring structures of $H^{*}(N^{n}(2);\field_{2})$ using the CW decomposition obtained by Kenso Fujii \cite{fuji1973onthekring}.
The above lemmas are shown in Section 3.
In Section 4, we show the main result.
The author express his gratitude to ...

\section{The cohomology ring of $N^n(m)$}

%$S^3/Q_8$は$N^n(m)$の特別な場合である$(n=0,m=2)$．
The spherical space form $S^3/Q_8$ is nothing but the manifold $N^n(m)$ introduced by K. Fujii \cite{fuji1973onthekring} in case when $n=0$ and $m=2$.
So we study the structure of the cohomology ring of $N^n(m)$ in a slightly general situation of $m=2$ and $n \ge 0$, 
%$N^n(m)$のCW複体としての構造やコホモロジー群は\cite{fuji1973onthekring}で議論されている．
using the cell-structure and the cohomology groups of general $N^n(m)$ obtained in \cite{fuji1973onthekring}.
%この節では\cite{fuji1973onthekring}の一部を紹介し，また$m=2$の場合のコホモロジー環を決定する．

\subsection{CW-structure of $N^n(m)$}

Let $\mathbb{H}$ be the skew-field of quaternion numbers, generated over $\real$ by $1$, $i$, $j$ and $k=ij$ with $i^{2}=j^{2}=-1$ and $ij=-ji$.
%一般四元数群$H_m$を$m\ge 2$に対して
We denote $H_m$ the generalised quaternion group for $m\ge 2$, which is defined by 

\[H_m = \langle x, y | x^{2^m} = y^4 = 1, y^2 = x^{2^{m-1}}, xyx = y \rangle.\]

%と定義する．
As a special case, we obtain $H_{2}=Q_{8}$ the quaternion group generated by $i$ and $j$ in $\quaternion$ which is the subgroup of $S^{3}$ the unit sphere of $\quaternion$.
%$H_2$は四元数群$Q_8$である．これは次のように$\mathbb{H}$の単位球面$S^3$の部分群とみなせる:
More generally, $H_{m}$ can be represented as the subgroup of $S^{3}$:

\[x = \exp (\pi i /2^{m-1}) \mbox{ and } y = j.\]

%$N^n(m)$ は$S^{4n+3}/H_m$ と定義される多様体である，
Using the above notions, Fujii introduced $N^{n}(m)$ in \cite{fuji1973onthekring} as the quotient manifold $S^{4n+3}/H_m$ from $S^{4n+3}= \{(q_1,\dots ,q_{n+1}) \in \quaternion^{n+1} \midvert \vert{q_1}\vert^{2} + \cdots + \vert{q_{n+1}}\vert^{2} = 1\}$ by the natural action of $H_{m}$ as follows:
%ここで$H_m$の作用は対角作用すなわち$S^{4n+3}$を$\mathbb{H}^{n+1}$内の単位球面とみなして
\[q(q_1,\dots ,q_{n+1}) = (qq_1,\dots ,qq_{n+1}),\]
%と作用する．
where $q \in S^{3}$ and $(q_1,\dots ,q_{n+1}) \in S^{4n+3} \subset \quaternion^{n+1}$.

\medskip

Since $S^{\infty} = \bigcup S^{4n+3}$ is contractible, we can easily see that $\bigcup_{i=0}^{\infty}N^{i}$ is nothing but $K(H_m,1)$ or, say, $BH_m$.

\medskip

%\cite{fuji1973onthekring}の中でFujiは$N^n(m)$の胞体分割を与えた．
In this paper, we use the following finite cell decomposition of the manifold $N^n(m)$ due to Fujii \cite{fuji1973onthekring}.
\begin{thm}[{\cite[Lemma 2.1]{fuji1973onthekring}}]\label{dec}
    The manifold $N^n(m)$ can be decomposed as the finite cell complex whose cells are given by $e^{4k+s}$ and $e^{4k+t}_i$ for $0\!\le\!k \!\le\!n$, \,$s\!=\!0,3$, \,$t\!=\!1,2$ and $i\!=\!1,2$, where $e^{4k+s}$ is a $4k{+}s$-cell and $ e^{4k+t}_i$ is a $4k{+}t$-cell.
	Moreover their boundary formulas are given as follows:
    \begin{align*}
    &\partial e^{4k} = 2^{m+1}e^{4k-1}_1,\ \ \partial e^{4k+1}_1 = \partial e^{4k+1}_2 = 0,\\
    &\partial e^{4k+2}_1 = 2^{m-1}e^{4k+1}_1 - 2e^{4k+1}_2,\ \ \partial e^{4k+2}_2 = 2e^{4k+1}_1,\ \ \partial e^{4k+3} = 0
    \end{align*}
\end{thm}

%この胞体分割の下で$N^{n-1}$は$N^n$の部分複体になっている．
Let us remark that the submanifold $N^{n-1}$ of $N^{n}$ is a subcomplex of $N_{n}$ with respect to the above cell decompostion.
%$\bigcup_{i=0}^{\infty}N^{i}$は$K(H_m,1)$すなわち$BH_m$である．

\subsection{The cohomology ring of $N^n(2)$}

%この補題から$N_n(m)$のホモロジー群を与える．
Firstly, we give a description of the homology groups of $N_{n}(m)$ by using the CW decomposition of $N_{n}(m)$ given in Theorem \ref{dec} due to Fujii \cite{fuji1973onthekring}.

\begin{prop}\label{prop:homologyandcohomology}
The homology and cohomology groups of $N^{n}(m)$ are given as follows, where $\Z_{d}$ denotes the cyclic group of order $d$ and $\F_{2}$ denotes the prime field of order $2$:
%~\\
%    \begin{enumerate}
%    \item The integral homology groups and cohomology groups of $N^{n}(m)$ are given as follows:
    \[
    H_k(N^n(m):\Z)=
    \begin{cases}
    \Z & k=0, 4n + 3,\\
    \Z_2 \oplus \Z_2 & k\equiv 1(4),  0<k<4n+3,\\
    \Z_{2^{m+1}} & k\equiv 3 (4),  0<k<4n+3,\\
    0 & otherwise.
    \end{cases}
    \]
    \[
    H^k(N^n(m):\Z)=
    \begin{cases}
    \Z & k = 0, 4n + 3,\\
    \Z_{2^{m + 1}} & k \equiv 0(4) , 0 <k<4n + 3,\\
    \Z_2\oplus \Z_2 & k \equiv 2(4), 0 <k< 4n + 3,\\
    0 & otherwise.
    \end{cases}
    \]
%    \item The cohomology groups of $N^n(m)$ and $N(m)$ with coefficients in $\F_2$ are given by
    \begin{align*}
    H^k(N^n(m):\F_2) &=
    \begin{cases}
    %\Z_2 & k=0, 4n + 3,\\
    \F_2 \oplus \F_2 & k\equiv 1,2(4),  0<k<4n+3,\\
    \F_2 & k\equiv 3,0 (4),  0\le k\le 4n+3,\\
    \end{cases}
    \\
    H^k(N(m): \F_2) &=
    \begin{cases}
    \F_2 \oplus \F_2 & k\equiv 1,2(4),\\
    \F_2 & k\equiv 3,0 (4).\\
    \end{cases}
    \end{align*}
%    \end{enumerate}
\end{prop}

%$H^{\ast}(N^n(2))$の積構造はファイブレーション
%\[S^3 \to X \to B\pi. \]
%に関するSerreスペクトル系列を考える，ここで$\pi = Q_8, X=S^3/\pi$である．
Secondly, we must determine the ring structure of $H^{\ast}(N^n(2))$ using Serre spectral sequence for the fibration 
\[S^3 \to X \to B\pi,\]
where $\pi = Q_8$ and $X = S^3/\pi$.

\begin{figure}[h!]
  \centering
  \includegraphics[width = 6cm]{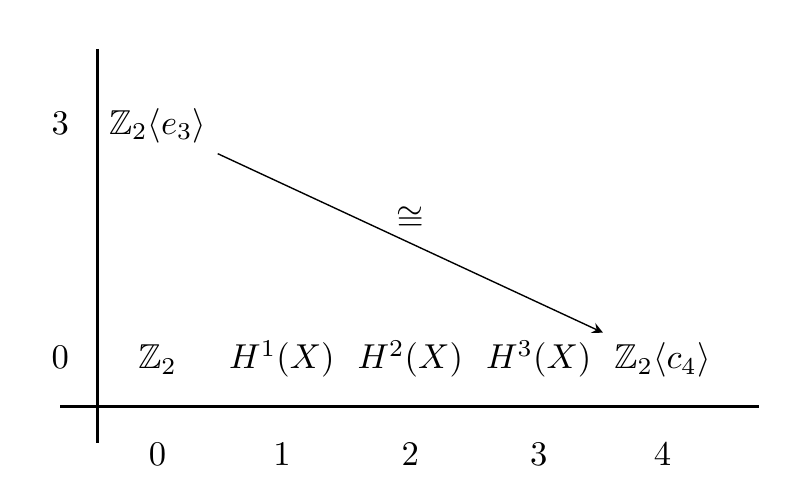}
  %[width = 6.5cm]
%   \caption{}
  %\label{}
\end{figure}
\noindent As is seen later in this section, this spectral sequence collapses at the fifth term.
\begin{prop}
    $H^{\ast}(B\pi;\F_2) \cong \F_2\langle x,y \rangle/(x^3,y^3,x^2{+}y^2{+}xy) \otimes \F_2\langle c_4\rangle$
\end{prop}
\begin{proof}
    %証明をいくつかのステップに分ける．\\
    The proof is divided into several steps.
    \par\medskip
    (step1)\, We show that Bockstein operation $\beta_2 : H^{1}(B\pi;\Z_2) \to H^{2}(B\pi;\Z_2)$ is an isomorphism.
    % bockstein 準同型の定義から，
    The extension of abelian groups 
    \[ 0 \to \Z \xrightarrow{2} \Z \xrightarrow{\pi} \Z_2\to 0,\]
    %の誘導する長完全列
    induces the following long exact sequence, since $H^{1}(B\pi;\Z)=H^{3}(B\pi;\Z)=0$.
    \[0 \to H^1(B\pi;\Z_2) \xrightarrow{\delta} H^2(B\pi;\Z) \xrightarrow{2_{\ast}} H^2(B\pi;\Z) \xrightarrow{\pi_{\ast}} H^2(B\pi;\Z_2) \to 0.\]
%    において$2_{\ast}$すなわち2倍写像は Proposition \ref{prop:homologyandcohomology} から零写像である．
Here, $2_{\ast}$ is the homomorphism induced from a twice map $2 : \Z \to \Z$ and is also a twice map.
Then by Proposition \ref{prop:homologyandcohomology}, we have $2_{\ast}=0$, and hence $\pi_{\ast}$ and $\delta$ are isomorphisms.
%    bockstein 準同型の定義から，$\beta_2 = \pi_{\ast} \circ \delta$は同型写像となりよって同型である．
Hence by the definition of Bockstein operation, we obtain that $\beta_2 = \pi_{\ast} \circ \delta$ is an isomorphism.
    \par\medskip
    (step2)\, $x^2 + y^2 + xy = x^3 = y^3= 0$.\\
%    $B\pi$をBar構成によって複体の構造を与える．
Let us first consider the cell structure given by the Bar resolution.
%    $x,y$の代表元をそれぞれ $x([x^ay^b]) = a, y([x^ay^b]) = b \mod 2$ とする．
Let $x$ and $y$ be cocycles given by the following formula.
\[x([x^ay^b]) = a,\quad y([x^ay^b]) = b \mod 2\]
Our gorl is to find $v \in C^{1}(B\pi)$ such that $\delta v = x^2 + xy + y^2$.
In fact, we can see that if we defined $v$ as follows:
\begin{align*}
    v(a) = \begin{cases}
    0 ~~a = e, x^3, y^3, x^3y,\\
    1 ~~a = x, y, xy, x^2.
    \end{cases}
\end{align*}
A equation $\beta(x^2) = \beta(\beta(x)) = 0$ indicates
\begin{align*}
    0 &= \beta(x^2 + y^2 + xy) = \beta(xy) = x^2y + xy^2\\
    0 &= x(x^2 + xy + y^2) = x^3 + x^2y + xy^2 = x^3\\
    0 &= (x^2 + xy + y^2)y = y^3
\end{align*}
Therefore, $x^2=xy^2$ is a generator of $H^3(B\pi;\mathbb{F}_2)$ and the proof ends with the standard arguments on the spectral sequence above.
\end{proof}
\begin{cor}
    $H^{\ast}(X;\F_2) = \Z_2\langle x, y\rangle/(x^3, y^3, x^2 + y^2 + xy)$
\end{cor}

\section{Topological complexity of $S^3/Q_8$}

In this section we give a proof for the equation $\tc{X=S^3/Q_8} = 6$. Let us consider the cohomology ring with coefficients in $\mathbb{F}_2$.

Firstly, we show the following:
\begin{prop}\label{caln2}
    $\cat(X) = 3$ and $5 \le \tc{X} \le 6$.
\end{prop}
\begin{proof}
    By proposition \ref{cupsum1}, \ref{cupsum2}, and \ref{cupsum3}, we obtain
    \[3\le \wgt{x^2y} \le \cat(X) \le \dim (X) = 3,\]
    and $5 \le \wgt{\bar{x}^3\bar{y}^2}\le \tc{X} \le 2\cat(X) = 6$, where $\bar{x}$ and $\bar{y}$ are zero divisors defined by 
    $\bar{x} = 1\otimes x + x\otimes 1$ and $\bar{y} = 1\otimes y + y \otimes 1$ respectively.
\end{proof}
From this proposition, our aim is to show that $\tc{X}$ is greater than $5$.
Therefore, it is sufficient to show that the homomorphism 
\[(e_5)^{\ast} : H^{\ast}(X\!\times\! X) \to H^{\ast}(P^{5}_{B}\Omega_{B}(X\!\times\! X))\]
is not injective.
To proceed further, let us consider the following diagram, where the vertical maps are induced from inclusion map $i: X\times X\to B\pi \times B\pi$.
\begin{figure}[h!]\label{diag1}
    \[\xymatrix{
P^{5}_{B}\Omega_{B}(X\!\times\! X)\ar[r] \ar[d] \ar@/^20pt/@{.>}[rr]^{e_{5}}& P^{\infty}_{B}\Omega_{B}(X\!\times\! X)\ar[r] & X\!\times \! X \ar[d]\\
P^{5}_{B}\Omega_{B}(B\pi\!\times\!B\pi)\ar[r] \ar@/^20pt/@{.>}[rr]^{e'_{5}} &P^{\infty}_{B}\Omega_{B}(B\pi\!\times\!B\pi)\ar[r]&B\pi\!\times\!B\pi
}\]\caption{}\end{figure}
Let $[z\otimes z]$ be a generator of $H^{6}(B\pi\!\times\!B\pi)$ such that $[i^{\ast}(z\!\otimes\!z)]$ is a generator of $H^{6}(X\!\times\!X)$.
If we can see that $(e'_5)^{\ast}(z\!\otimes\!z)$ vanishes in $H^{6}(P^{5}_{B}\Omega_{B}(X\!\times\! X))$, the homomorphism $(e_5)^{\ast}$ is not injective.
The following lemma is well known (cf. M. C. Crabb and W. A. Sutherland \cite[prop 3.1]{crabb2000counting}), when $m=\infty$.
\begin{lem}
    The fibrewise projective $m$-space $P^{m}_{B}\Omega_{B}(B\pi\!\times\!B\pi)$ is fibrewise homotopic to $P^{m}_{B}\widehat{K} = S^{\infty}\times_{\text{ad}}P^m\pi$.
    In this case, we assume that $\widehat{K} = S^{\infty} \times_{\mbox{ad}}\pi$, $P^m\pi$ is the $m$-skelton of $B\pi$ with the cell structure given by the Bar resolution,
    and the adjoint action of $\pi$ given as follows:
    \[ h[g_1|g_2|\cdots | g_m] = [hg_1h^{-1}|hg_2h^{-1}|\cdots | hg_mh^{-1}],~ h \in \pi,~ [g_1|g_2|\cdots | g_m] \in P^{m}\pi .\]   
\end{lem}
\begin{proof}
    Let $\mathcal{P}(B\pi)  = \{ u : I \to S^{\infty} | p(u(0)) = p(u(1)) \}$, where $p : S^{\infty} \twoheadrightarrow B\pi$.
    Since $S^{\infty}$ is contractible, there is a homotopy $H : S^{\infty} \times I \to S^{\infty}$ from identity map to constant map, 
    and we describe adjoint map $\text{ad} H : S^{\infty} \to \text{Map}(I, S^{\infty})$ as $\mu$. 
    We define $f : \mathcal{P}(B\pi) \to \widehat{K}$ as $f(u) := [u_0 , [(p_{\ast}\mu_{u_0})^{-1} \cdot (p_{\ast} u) \cdot (p_{\ast} \mu_{u_1})]]$. 
    Moreover, when $p_{\ast}u = p_{\ast}v$,
\begin{align*}
     f(u) &= [u(0) , [(p_{\ast}\mu_{u(0)})^{-1} \cdot (p_{\ast} u) \cdot (p_{\ast} \mu_{u(1)})]]\\
     &= [v(0)\cdot {}^\exists h^{-1} , [(p_{\ast}\mu_{v(0)})^{-1}(p_{\ast}h)^{-1} \cdot (p_{\ast} v) \cdot (p_{\ast}h)(p_{\ast} \mu_{v(1)})]]\\
     &= [v(0) , [(p_{\ast}\mu_{v(0)})^{-1}(p_{\ast}h)^{-1}(p_{\ast}h) \cdot (p_{\ast} v) \cdot (p_{\ast}h)^{-1}(p_{\ast}h)(p_{\ast} \mu_{v(1)})]]\\
     &= [v(0) , [(p_{\ast}\mu_{v(0)})^{-1} \cdot (p_{\ast} v) \cdot (p_{\ast} \mu_{v(1)})]] = f(v),
\end{align*}
    so $f$ induces $\tilde{f} : \Omega_{B}(B\pi\!\times\!B\pi) \to \widehat{K}$, and the restriction of $\tilde{f}$ to each fibre is homotopy equivalence since $B\pi$ is a $K(\pi, 1)$. 
    Then by \cite{james1995introduction}, $\tilde{f}$ is a fibrewise homotopy equivalence. 
    In addition, $\tilde{f}$ is a fibrewise $A_{\infty}$-map since each fibre of $\widehat{K} \to B\pi$ is a discrete set. 
    Therefore, $P^{m}_{B}\Omega_{B}(B\pi\times B\pi) \simeq_{B} P^{m}_{B}\widehat{K}$.
\end{proof}

The cell structure of $P^{m}_{B}\widehat{K}$ is discribed by product cells of $S^{\infty}$ and $P^{m}\pi$.
\[ P^{m}_{B}\widehat{K}  = \bigcup_{0\le n \le m}\bigcup_{\omega \in \pi^n} \bigcup_{\sigma \in\{ \text{cells of } B\pi\} = \Lambda } [\sigma|\{\omega\}]. \]
Next, we describe the boundary formula up to dimension $3$.
By Lemma 2.1, the $3$-skelton of $S^{\infty}$ can be visualised as follows.
\begin{figure}[h!]
    \includegraphics[width = 6cm]{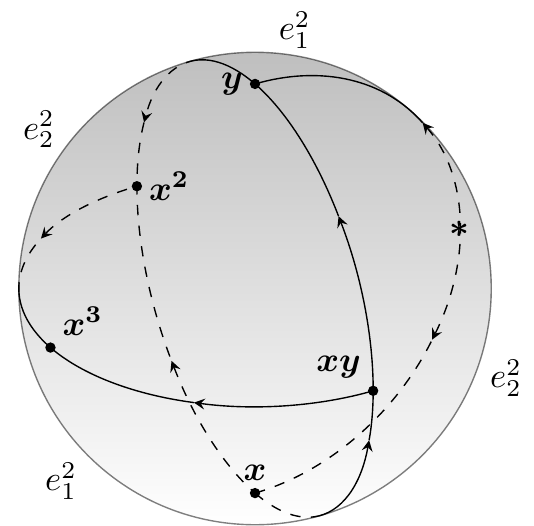}
\end{figure}
The boundary of a $k{+}n$-cell $[\eta|{\omega}]$ is a sum of $[\eta|\{\partial_i\omega\}]~(0\le i \le n)$, $[\eta'|\{h\omega h^{-1}\}]~(\eta'\in \Lambda_{k-1}, h \in \Lambda_1)$.
\begin{prop}
    The modulo 2 boundary of product cells \[ [e^{3} | \{ \omega \}],  [e^{2}_{1} | \{ \omega \}],  [e^{2}_{2} | \{ \omega \}],  [e^{1}_{1} | \{ \omega \}],  [e^{1}_{2} | \{ \omega \}],~ \omega \in \pi^n\] are given as follows.
    \begin{align*}
    \partial [e^3 | \{ \omega \}] &=[e^{2}_{1} | \{ \omega \}] + [e^{2}_{2} | \{ \omega \}] + [e^{2}_{1} | \{ x\bar{y}\omega y\bar{x} \}] + [e^{2}_{1} | \{ \bar{x} \omega x \}] + [e^3 | \{\partial\omega\}], \\
    \partial [e^{2}_{1} | \{ \omega \}] &=[e^{1}_{1} | \{ \omega \}] + [e^{1}_{2} | \{ \omega \}] + [e^{1}_{1} | \{\bar{x} \omega x \}] + [e^{1}_{2} | \{\bar{y} \omega y \}] + [e^{2}_{1} | \{\partial\omega\}], \\
    \partial [e^{2}_{2} | \{ \omega \}] &=[e^{1}_{1} | \{ \omega \}] + [e^{1}_{2} | \{ \omega \}] + [e^{1}_{1} | \{\bar{y}\bar{x} \omega xy \}] + [e^{1}_{2} | \{\bar{x} \omega x \}] + [e^{2}_{2} | \{\partial\omega\}], \\
    \partial [e^{1}_{1} | \{ \omega \}] &=[\ast | \{ \omega \}] + [\ast | \{ \bar{x}\omega x \}] + [e^{1}_{1} | \{\partial\omega\}], \\
    \partial [e^{1}_{2} | \{ \omega \}] &=[\ast | \{ \omega \}] + [\ast | \{ \bar{y}\omega y \}] + [e^{1}_{2}| \{\partial\omega\}],
    \end{align*}
    where $\displaystyle[\sigma | \{\partial\omega\}] = \overset{n}{\underset{i = 0}{\textstyle\sum}} (-1)^{i}[\sigma | \{ \partial_{i} \omega \}]$.
\end{prop}
\begin{proof}
    The $3$-skelton of $X$ can be described as follows: $\partial [e^3 | \{ \omega \}]$ is the union of cells
     $[e^{2}_{1} | \{ \omega \}],
    [e^{2}_{2} | \{ \omega \}],
    [e^{2}_{1}y\bar{x} | \{\omega \}] = [e^{2}_{1} | \{ x\bar{y}\omega y\bar{x} \}],
    [e^{2}_{1}x | \{ \omega \}] = [e^{2}_{1} | \{ \bar{x} \omega x \}]$, and
    $[e^3 | \{\partial\omega\}]$ by the previously given CW structure of $X$.
%    Similarly, we obtain 
    The boundary formulas of the product cells $[e^{2}_{1} | \{ \omega \}]$, $[e^{2}_{2} | \{ \omega \}]$, $ [e^{1}_{1} | \{ \omega \}]$, $[e^{1}_{2} | \{ \omega \}]$ can be described similarly.
    For example, the boundary of $e_{1}^{2}$ is the union of cells $[e^{1}_{1}|\{\omega\}], [e^{1}_{2}|\{\omega\}], [e^{1}_{2}x|\{\omega\}]$, and $[e^{1}_{1}xy|\{\omega\}]$
    therefore, $\partial [e^{2}_{2}|\{\omega\}] = [e^{1}_{1}|\{\omega\}]+ [e^{1}_{2}|\{\omega\}]+ [e^{1}_{2}|\{\bar{x}\omega x\}]+ [e^{1}_{1}|\{\bar{x}\bar{y}\omega xy\}]+ [e^{2}_{2}|\{\partial\omega\}]$.
    Therefore we obtain the proposition.
    % {\color{red} Why don't you describe more precisely?}
\end{proof}

\newcommand{\pullbackmark}[2]{\save ;p+<1.4pc,0pc>:(0,-1)::%記号の大きさを指定
(#1) *{\phantom{Z}} %
;p+(#2)-(0,0) **@{-}%
;p-(#1)+(0,0) *{\phantom{Z}} **@{-} \restore}
We can update Diag.1 using lemma 3.2.
\begin{figure}[h!]
\[\xymatrix{
&\ar[dl] P^{5}_{B}\Omega_{B}(X\!\times\!X)\ar[r] \ar[d] \ar@/^20pt/@{.>}[rr]^{e_{5}}& P^{\infty}_{B}\Omega_{B}(X\!\times\!X)\ar[r] & X\!\times \! X \ar[dd]\\
X\ar[d] & \ar[l] S^{3}\times_{\text{ad}} P^{5}\pi \ar[d]^{\iota} \pullbackmark{-1,0}{0,1} & & \\
B\pi&\ar[l] S^{\infty}\times_{\text{ad}}P^{5}\pi \ar[r] \ar@/^20pt/@{.>}[rr]^{e'_{5}} & \ar[r]^-{\simeq} S^{\infty}\times_{\text{ad}}P^{\infty}\pi &B\pi\!\times\!B\pi
}\]\caption{}
\end{figure}
Since We are replacing the space with another space with homotopy equivalence, the above figure is not commutative.
However, it does not cause problems when discussing the injectivity of the homomorphism $(e_5)^{\ast}$.
In particular, the Diag.2 become commutative because we concider cohomology with $\mathbb{F}_2$ coefficient.

We can discribe $(\iota)^{\ast}(e'_5)^{\ast}(z\!\otimes\!z)$ explicitly:
\[ (\iota)^{\ast}(e'_{5})^{\ast}(z\otimes z)[\sigma| \{ h_1 | \dots | h_m\}] = \begin{cases} x^2y[h_1 | h_2 | h_3] & \text{if $\sigma = e^3$, $m= 3$},\\
0 &\text{otherwise}. \end{cases} \]

Crealy, the boundary homomorphism $\delta: C^{5}(S^3\times_{\mbox{ad}}P^5\pi) \to C^{6}(S^3\times_{\mbox{ad}}P^5\pi)$ is a liner map between vector spaces $C^{5}(S^3\times_{\mbox{ad}}P^5\pi)$
and $C^{6}(S^3\times_{\mbox{ad}}P^5\pi)$ whose basises are the set of $5$ and $6$-cells $[\sigma | \omega], \sigma \in S^3/Q_8, \omega \in P^{\ast}\pi$ respectively.
Therefore, in order to find a $u \in C^{5}(S^3\times_{\mbox{ad}}P^5\pi)$ such that $\delta u = (\iota)^{\ast}(e'_{5})^{\ast}(z\otimes z)$, 
we just calculate the rank of $\delta$.
However, this calculation is very hard even if we use a computer, so we suppose that one element satisfying $\delta u = (\iota)^{\ast}(e'_{5})^{\ast}(z\otimes z)$ can discribed as cup product of $v \in C^2(S^3\times_{\mbox{ad}}P^5\pi)$
\[ v[\sigma| \{ h_1 | \dots | h_m\}] = \begin{cases} xy[h_1 | h_2] & \text{if $\sigma = e^3$, $m= 2$},\\
0 &\text{otherwise} \end{cases} \]
and $u'\in C^2(S^3\times_{\mbox{ad}}P^5\pi)$ such that
\begin{align}\label{A}
(\delta u')[\sigma| \{ h_1 | \dots | h_m\}] = \begin{cases} x[h] & \text{if $\sigma = e^3$, $m= 1$},\\ 0 &\text{otherwise}. \end{cases} \tag{A}
\end{align}
In this case, we can solve the simultaneous equations $\delta u' = (\iota)^{\ast}(e'_{5})^{\ast}x$ using a computer and in fact, such $u'$ exists.
\begin{thm}
    $\tc{X} = \wgt[B]{X} = \wgt[B]{z\otimes z} = 6$.
\end{thm}
\begin{proof}
    By the above argument, we have \[(\iota)^{\ast}(e_{5})^{\ast}(z\otimes z) = [\delta u] = [\delta (u' \smile v)] = [0]\]
    and $z\otimes z \neq 0$. So it implies $\wgt[B]{X} \ge 6$. On the other hand, Theorem \ref{caln2} we have $\tc{X} \le 2\cat(X) \le 6$ therefore, $6 \le \wgt[B]{X} \le \tc{X} \le 2\cat(X) = 6$.
\end{proof}
\section*{Appendix A}
I write the algorithm that I made to see the existense of such $u'$ satisfying the equation $(A)$ in section 3.

\begin{algorithm}[H]
    \SetAlgoLined
    \KwData{$G$ : Cayley table of $Q_8$, $\Lambda$ : the set of cells of $S^3/Q_8$, the modulo 2 boundary formula of product cells and $\omega$ : cohomology class in $H^{4}(S^3\times_{\text{ad}}P^5\pi)$}
    \KwResult{the augmented matrix $A$ corresponding the equation $(A)$}
    $C_4$ : The set of 4-cells of $S^3\times_{\text{ad}}P^5\pi$\\
    $V = \{\}$ /* The set of variables */\\
    \For{ $x \in C_4$ }{
      $\partial x$ : the set of cells come out from the modulo 2 boundary of $x$\\
      $V = V + \partial x$
      }
    $F = []$ : /* the list of equations $\{ f_i = b_i\}_i$ */\\
    $column$ = [the list of variables that is in $V$] + ['$answer$']\\
    \For{$x \in C_4$}{
       $(f_x | b_x) := $ the tuple with column '$column$' \\
       \eIf{$\omega(x) = 1$}{$(f_x | b_x)_{answer} = 1$}{$(f_x | b_x)_{answer} = 0$}
       \For{$y \in V$}{
        \eIf{$y \in \partial x$}{
           $(f_x | b_x)_y = 1$
        }{
           $(f_x | b_x)_y = 0$
        }
       }
       /* the tuple $(f_x | b_x)$ consists of 0 or 1 */\\
       add ${}^t(f_x | b_x)$ to $F$}
       $A = {}^tF = {}^t[{}^t(f_{x_1} | b_{x_1}),\dots,{}^t(f_{x_n} | b_{x_n})]_{x_i \in C_4}$\\
    
    \caption{algorithms to create augmented matrix $A$}
  \end{algorithm}
\section*{Acknowledgement}
I would like to thanks to Norio Iwase. This paper would not exist without his detailed and valuable advice.
%%%%%%%%%参考文献%%%%%%%%%%%%%
%
%	References
%
\bibliographystyle{alpha}
\bibliography{sankou}
%%%%%%%%%%%%%%%%%%%%%%%%%%%%%
\end{document}